\documentclass[a4paper,12pt]{article}

\usepackage{natbib}             
\usepackage{graphicx}           
\usepackage{amsfonts}
\usepackage{amsmath}

\def\R{\mathbb{R}}

\begin{document}

\title{\bf Measures and LMI for space launcher robust control validation\footnote{This work was performed
in the frame of the SAFE-V (Space Application Flight control Enhancement of Validation Framework)
ESA (European Science Agency) TRP (Technology Development Project) under contract number 4000102288.
The research of the first author was also supported by research project 103/10/0628 of the Grant Agency
of the Czech Republic.}}

\author{Didier Henrion$^{1,2,3}$, Martine Ganet-Schoeller$^4$, Samir Bennani$^5$}

\footnotetext[1]{CNRS; LAAS; 7 avenue du colonel Roche, F-31400 Toulouse; France. {\tt henrion@laas.fr}}
\footnotetext[2]{Universit\'e de Toulouse;LAAS; F-31400 Toulouse; France.}
\footnotetext[3]{Faculty of Electrical Engineering, Czech Technical University in Prague,
Technick\'a 2, CZ-16626 Prague, Czech Republic.}
\footnotetext[4]{EADS ASTRIUM Space Transportation; 66 route de Verneuil, F-78133~Les Mureaux; France.}
\footnotetext[5]{ESA/ESTEC (TEC-ECN); Guidance, Navigation \& Control Group; Keplerlaan 1, NL-2201~AZ Noordwijk; The~Netherlands.}

\maketitle

\begin{abstract}
We describe a new temporal verification framework for safety and robustness
analysis of nonlinear control laws, our target application being a space launcher vehicle.
Robustness analysis, formulated as a nonconvex nonlinear optimization problem
on admissible trajectories corresponding to piecewise polynomial dynamics,
is relaxed into a convex linear programming problem on measures.
This infinite-dimensional problem is then formulated as a generalized moment problem,
which allows for a numerical solution via a hierarchy of linear matrix inequality
relaxations solved by semidefinite programming. The approach is illustrated
on space launcher vehicle benchmark problems, in the presence of closed-loop
nonlinearities (saturations and dead-zones) and axis coupling.
\end{abstract}

\begin{center}
\end{center}

\section{Introduction}

This work is carried out within the scope of the project SAFE-V (Space Application
Flight control Enhancement of Validation Framework). The objective of this project
is to analyse, develop and demonstrate effective design, verification and validation
strategies and metrics for advanced guidance navigation and control systems.
These new strategies should be implemented in an incremental manner in the traditional
validation framework in order to be applicable both for current launchers validation
and for future launcher and re-entry vehicle validation.

Control design performance and compliance with system specifications should be demonstrated
by a suitable combination of stability analysis and simulation. This validation is made
on a wide domain of variation of influent parameters; it includes worst case validation
and Monte Carlo simulation for statistic requirements, see e.g. \cite{RON99}
for Ariane 5 validation.

Generally speaking, the Monte Carlo method is a numerical integration method using sampling, which can be used,
for example, to determine the quantities of interest of a variable of interest
for a random input variable such as mean or standard deviation or probability density function.
One of the most commonly used methods inside ASTRIUM ST for the confidence
interval on a quantile estimation is Wilks formula, see \cite{OPE08}. This formula of the
confidence interval depending on the probability level and the number of samples can be used,
for Monte Carlo only, in two ways: to determine the probability and confidence level for the
values of samples chosen by the user, or in reverse to determine the number of
simulations to be carried out for the values probability and confidence level chosen by the user.
Monte Carlo simulations are widely used in traditional validation process and they are very robust
because they do not require smoothing assumptions on the model neither importance factor
assumptions on the parameters. The random events are simulated as they would occur naturally.
The generic Monte Carlo method is however very time consuming,  because a great number
of simulations to demonstrate extreme level of probability (in general much more than
$1/p$ simulations to estimate the probability of an event of probability $p$).
The number of simulations required may be critical for long duration mission simulation
(like launcher flight). In the context of robust control, see \cite{tempo05}
for related probabilistic approaches. Note that these approaches may miss worst
case behavior, especially when the number of parameters is large.

In contrast, worst case validation techniques are aimed at identifying worst
case instances and behavior. These techniques include Lyapunov and mu-analysis approaches.
Among them, the most relevant for space vehicle applications are:
\begin{itemize} 
\item results by the GARTEUR group, see \cite{FIE02}, and the COFCLUO group, see \cite{COF09},
for recent aeronautic applications and new developments;
\item results by ASTRIUM ST in the frame of ESA ITT for Robust LPV for launcher application, see
\cite{TN069,TN079,GAN11}, and for robust stability analysis for ATV industrial validation, see \cite{GAN09};
\item results by NGC on rendezvous validation process in the frame of the ESA VVAF project, see
\cite{disotto,disotto2,KRO10-1,paulino},
and on re-entry validation, see \cite{KRO04}.
\end{itemize}
Many other non industrial applications are available in the literature. Two of them could be
considered as major for our applications: a PhD thesis on missile application by \cite{ADO04}
and the IQC toolbox developed by \cite{SCH06} with Delft University that was applied
for HL-20 re-entry vehicle analysis in \cite{VEN09}, under ESA contract. 

In this paper we propose an original approach to validation of control laws based on
measures and convex optimization over linear matrix inequalities (LMIs). The approach
can be seen as a blend between worst case validation techniques and statistical simulation
techniques. On the one hand, we optimize over worst case admissible trajectories
to validate a system property. On the other hand, we propagate along system trajectories
statistical information (probability measures) instead of deterministic initial conditions.
Moreover, our approach is primal in the sense that we optimize directly over systems
trajectories, we do not seek a dual Lyapunov certificate.

First we rephrase our validation problem as a
robustness analysis problem, and then as a nonconvex nonlinear optimization
problems over admissible trajectories. This is the approach followed e.g. in
\cite{prajna07} where the authors verify or prove temporal properties such as
safety (all trajectories starting from a set of initial conditions
never reach a set of bad states), avoidance (at least one trajectory starting from
a set of initial conditions will never reach a set of bad states), eventuality
(all trajectories starting from a set of initial conditions will reach a set of
good states in finite time) and reachability (at least one trajectory starting from a set of
initial conditions will reach a set of good states in finite time).

Following \cite{lasserre08}, we formulate our nonconvex nonlinear trajectory optimization problem 
as a linear programming (LP) problem on measures, with at most three unknowns: the initial
measure (modeling initial conditions), the occupation measure (encoding system trajectories) and
the terminal measure (modeling terminal conditions); in some cases
a subset of these decision variables may be given. The final time can be given or free,
finite or infinite. If all the data are polynomials
(performance measure, dynamics, constraints), then we formulate the infinite-dimensional LP
on measures as a generalized moment problem (GMP), and we approach it via
an asymptotically converging hierarchy of finite-dimensional convex LMI relaxations.
These LMI relaxations are then modeled with our specialized
software GloptiPoly 3, solved with a general-purpose
semidefinite programming (SDP) solver, and system properties are validated from
the numerical solutions.

The models we can deal with should be piecewise polynomial (or rational) in time and space.
We consider a compact region of the state-space over which the trajectories are optimized.
The system dynamics are defined locally on explicitly given basic semialgebraic sets
(intersections of polynomial sublevel sets) by polynomial vector fields. The objective
function (stability or performance measure) consists of a polynomial terminal term
and of the time integral of a polynomial integrand.

Our contribution is to use the approach described in \cite{lasserre08} jointly
with the corresponding software tool GloptiPoly 3 by \cite{henrion09} to address
temporal validation/verification of control systems properties, in the
applied context of space launcher applications. The use of LMI and measures
was already investigated in \cite{prajna07} for building Lyapunov barrier
certificates, and based on a dual to Lyapunov's theorem described
in \cite{rantzer01}. Our approach is similar, in the sense that optimization over systems
trajectories is formulated as an LP in the infinite-dimensional space
of measures. This LP problem is then approached as a generalized
moment problem via a hierarchy of LMI relaxations.

Historically, the idea of reformulating nonconvex nonlinear ordinary differential equations (ODE) into 
convex LP, and especially linear partial differential equations (PDE)
in the space of probability measures, is not new. 
It was Joseph Liouville in 1838 who first introduced the linear PDE involving the Jacobian
of the transformation exerted by the solution of an ODE on its initial condition \cite{liouville},
see \cite{ehrendorfer94,ehrendorfer02} for a survey on the Liouville PDE with applications
in meteorology. The idea was then largely expanded in Henri Poincar\'e's work on dynamical systems
at the end of the 19th century, see in particular 
\cite[Chapitre XII (Invariants int\'egraux)]{poincare99}. This work was pursued in the 20th
century in \cite{kryloff37}, \cite[Chapter VI (Systems with an integral invariant)]{nemystkii47}
and more recently in the context of optimal transport by e.g. \cite{rachev98}, \cite{villani03}
or \cite{ambrosio}.
Poincar\'e himself in \cite[Section IV]{poincare08} mentions the potential of formulating
nonlinear ODEs as linear PDEs, and this programme has been carried out to
some extent by \cite{carleman32}, see also \cite{lasota85}, \cite{kowalski91},
\cite{hernandez03} and more recently \cite{barkley06}, \cite{vaidya08}, \cite{gaitsgory09}.
For recent studies of the Liouville equation in optimal control see e.g. \cite{kwee}
and \cite{brockett}, and for applications in uncertainty propagation and control law
validation see e.g. \cite{mellodge} and \cite{halder}.

\section{Piecewise polynomial dynamic optimization}

Consider the following dynamic optimization problem with piecewise polynomial differential constraints
\begin{equation}\label{ode}
\begin{array}{lll}
J \:= & \inf_{x(t)} & h_T(T, x(T)) + \int_0^T h(t,x(t))dt \\
& \mathrm{s.t.} & \dot{x}(t) = f_k(t,x(t)), \:\: x(t) \in {X}_k, \:\: k=1,2,\ldots,N \\
& & x(0) \in {X}_0, \:\: x(T) \in {X}_T, \:\: t \in [0,T]
\end{array}
\end{equation}
with given polynomial dynamics $f_k \in \R[t,x]$ and costs
$h, h_T \in \R[t,x]$, and state trajectory $x(t)$
constrained in closed basic semialgebraic sets
\[
{X}_k = \{x \in \R^n \: : \: g_{kj}(t,x) \geq 0, \:
j=1,2,\ldots,N_k\}
\]
for given polynomials $g_{kj} \in \R[t,x]$.
We assume that the state-space partitioning sets, or cells ${X}_k$, are disjoint,
i.e. all their respective intersections have zero Lebesgue measure in $\R^n$, and they all belong to
a given compact semialgebraic set, e.g.
\[
{X} = \{x \in \R^n \: :\: \|x\|^2_2 \leq M\}
\]
for a sufficiently large constant $M>0$.
Finally, initial and terminal states are constrained in semialgebraic ets
\[
{X}_0 = \{x \in \R^n \: : \: g_{0j}(t,x) \geq 0, \:
j=1,2,\ldots,N_0\} \subset {X}
\]
and
\[
{X}_T = \{x \in \R^n \: : \: g_{Tj}(t,x) \geq 0, \:
j=1,2,\ldots,N_T\} \subset {X}
\]
for given polynomials $g_{0j}, g_{Tj} \in \R[t,x]$.

It is assumed that optimization problem (\ref{ode}) arises from validation
of a controlled (closed-loop) system behavior. In this problem, the infimum
is sought over absolutely continuous state trajectories
$x(t)$, but the same methodology can be extended to (possibly
discontinuous) trajectories of bounded variation (e.g. in the presence
of impulsive controls or state jumps), see \cite{claeys11}.

The final time $T$ is either given, or free, in which case it becomes
a decision variable, jointly with $x(t)$. Similarly, the initial and
terminal constraint sets ${X}_0$ and ${X}_T$ are either given, or free, in which
case they also become optimized decision variables.

\section{Linear programming on measures}

In this section, we formulate optimization problem (\ref{ode}),
which is nonlinear and nonconvex in state trajectory $x$,
into a convex optimization problem linear in $\mu$, the occupation
measure of trajectory $x$. This reformulation is classical
in the modern theory of calculus of variations and optimal control
and it can be traced back to \cite{young37}, see also \cite{ghouila67},
\cite{young69} and \cite{gamkrelidze75}, amongst many others.

\subsection{Occupation measure}

To understand the basic idea behind the transformation and the
elementary concept of occupation measure, it is better to deal with the
single nonlinear ODE
\begin{equation}\label{sode}
\dot{x}_t = f(x_t)
\end{equation}
where $x_t$ is a shortcut for $x(t)$, a time-dependent state vector
of $\R^n$ and $f : \R^n \to \R^n$ is a uniformly Lipschitz map. It follows that
the Cauchy problem for ODE (\ref{sode}) has a unique solution $x_t$
for any given initial condition $x_0 \in \R^n$.

Now think of initial condition $x_0$ as a random variable of $\R^n$, or more
abstractly as a nonnegative probability measure $\mu_0$ with support ${X}_0 \subset \R^n$,
that is a map from the sigma-algebra of subsets of ${X}_0$ to the interval $[0,1] \subset \R$
such that $\mu_0({X}_0)=1$. For example, the expected value of $x_0$ is the vector
$E[x_0] = \int_{{X}_0} x \mu_0(dx)$.

Now solve ODE (\ref{sode}) for a trajectory, or flow $x_t$, given this
random initial condition. At each time $t$, state $x_t$ can also be interpreted
as a random variable, i.e. a probability measure that we denote by $\mu_t(dx)$. We say that the measure
is transported by the flow of the ODE. The one-dimensional family, or path of measures $\mu_t$ satisfies a
PDE
\begin{equation}\label{pde}
\frac{\partial \mu_t}{\partial t} + \mathrm{div}(f \mu_t) = 0
\end{equation}
which turns out to be linear in the space of probability measures.
This PDE is usually called Liouville's equation.
As explained e.g. in \cite[Theorem 5.34]{villani03}, nonlinear ODE (\ref{ode})
follows by applying Cauchy's method of characteristics to linear transport PDE (\ref{pde}),
see also \cite[Section 3.2]{evans10} for a tutorial exposition.
In equation (\ref{pde}), div models the divergence operator, i.e.
\[
\mathrm{div}(v) = \sum_{i=1}^n \frac{\partial v_i}{\partial x_i}
\]
for every smooth function $v$. Its action on measures should be understood
in the weak, or distributional sense, i.e.
\[
\int v\:\mathrm{div}(\nu) = -\int Dv \cdot d\nu
\]
where $v$ is a smooth test function, $\nu$ is a vector-valued measure
and $D$ stands for the gradient operator.
Given a subset ${\mathcal T}\times{\mathcal X}$ in the sigma-algebra of subsets of
$[0,T]\times X$, we define the occupation measure
\[
\mu({\mathcal T}\times{\mathcal X}) = \int_{\mathcal T} \mu_t({\mathcal X})dt
\]
which encodes the time-space trajectories $x_t$, in the sense
that $\mu([0,T]\times{\mathcal X})$ is the total time spent by
trajectory $x_t$ in a subset ${\mathcal X}\subset X$ of the state space.
In its integral form, transport PDE (\ref{pde}) becomes
\begin{equation}\label{ipde}
\mathrm{div}(f \mu) = \mu_0 - \mu_T
\end{equation}
where $\mu_T$ is the terminal probability measure with support ${X}_T \subset \R^n$.
PDE (\ref{ipde}) can equivalently be formulated as
\begin{equation}\label{tpde}
\int_{X} Dv \cdot f d\mu = \int_{{X}_T} v d\mu_T - \int_{{X}_0} v d\mu_0 
\end{equation}
for all smooth functions $v$ compactly supported on ${X}$.
Problem (\ref{tpde}) is an infinite-dimensional linear system
of equations linking occupation measure $\mu$, initial measure $\mu_0$
and terminal measure $\mu_T$, consistently with ODE (\ref{sode}).

\subsection{Finite terminal time}

If we apply these ideas to problem (\ref{ode}), encoding the state trajectory $x(t)$
in an occupation measure $\mu$ supported on $X$, we come up with
an infinite-dimensional LP problem
\begin{equation}\label{lp}
\begin{array}{lll}
J_{\infty} \:= & \inf_{\mu} & \int h_T(T,x) d\mu_T(x) + \sum_k \int h(t,x) d\mu_k(t,x) \\
& \mathrm{s.t.} & \sum_k \int \frac{\partial v(t,x)}{\partial t} \:d\mu_k(t,x) \:\:+ \\
&& \sum_k \int Dv(t,x) \cdot f_k(t,x) \:d\mu_k(t,x) \:\:= \\
&& \int v(T,x) \:d\mu_T(x) - \int v(0,x) \:d\mu_0(x) 
\end{array}
\end{equation}
for all smooth test functions $v$,
where each occupation measure $\mu_k$ is supported on set ${X}_k$
and the global occupation measure is
\[
\mu = \sum_k \mu_k
\]
with normalization constraint
\begin{equation}\label{norm}
\mu([0,T]\times{X}) = T
\end{equation}
such that $T$ is a finite terminal time.

In problem (\ref{lp}), final time $T$, initial measure $\mu_0$ and terminal
measure $\mu_T$ may be given, or unknown, depending on the original optimization
problem. 

\subsection{Generalized moment problem}

More concisely, LP problem (\ref{lp}) can be formulated as follows:
\begin{equation}\label{measlp}
\begin{array}{ll}
\inf_{\mu} & \sum_k \int_{{X}_k} c_k\:d\mu_k \\
\mathrm{s.t.} & \sum_k \int_{{X}_k} a_{ki}\:d\mu_k = b_i, \quad\forall i
\end{array}
\end{equation}
where the unknowns are a finite set of nonnegative measures $\mu_k$,
with respective compact semialgebraic supports
\begin{equation}\label{sets}
{X}_k = \{x \: :\: g_{kj}(x) \geq 0, \:\:\forall j\}.
\end{equation}
Note that here we have incorporated time variable $t$ into vector $x$,
for notational conciseness.
If all the data in problem (\ref{lp}) are polynomials, and if we
generate test functions $v(x)$ using a polynomial basis (e.g. monomials,
which are dense in the set of continuous functions with compact support),
all the coefficients $a(x)$, $b(x)$, $c(x)$ are polynomials, and there
is an infinite but countable number of linear constraints indexed by $i$.

We will then manipulate each measure $\mu_k$ via its moments
\begin{equation}\label{mom}
y_{k\alpha} = \int_{{X}_k} x^{\alpha} d\mu_k(x), \quad \forall \alpha
\end{equation}
gathered into an infinite-dimensional sequence $y_k$ indexed by a vector
of integers $\alpha$,
where we use the multi-index notation $x^{\alpha} = x_1^{\alpha_1} x_2^{\alpha_2} \ldots $
LP measure problem (\ref{measlp}) becomes an LP moment problem, or GMP, see
\cite{lasserre09}:
\[
\begin{array}{ll}
\inf_y & \sum_k \sum_{\alpha} c_{k\alpha} y_{k\alpha} \\
\mathrm{s.t.} & \sum_k \sum_{\alpha} a_{ki\alpha} y_{k\alpha} = b_i, \quad \forall i
\end{array}
\]
provided we can handle the representation condition (\ref{mom})
which links a measure with its moments. It turns out \cite[Chapter 3]{lasserre09}
that if sets ${X}_k$ are compact semialgebraic as in (\ref{sets}),
we can use results of functional analysis and real algebraic geometry
to design a hierarchy of LMI relaxations which is asymptotically
equivalent to the generalized moment problem. In each LMI relaxation
\begin{equation}\label{lmi}
\begin{array}{lll}
J_d \:= & \inf_y & c^T y \\
& \mathrm{s.t.} & Ay = b \\
& & M_d(y) \succeq 0 \\
& & M_d(g_{kj}, y) \succeq 0, \quad\forall j,k
\end{array}
\end{equation}
we truncate the moment sequence to a finite number of its moments.
If the highest moment power is $2d$, we call LMI relaxation of order $d$
the resulting finite-dimensional truncation. Matrix $M_d(y)$ is symmetric and
linear in $y$, it is called the moment matrix, and it must be
positive semidefinite for (\ref{mom}) to hold. Symmetric matrices $M_d(g_{kj}, y)$
are also linear in $y$ and they are called localizing matrices. They ensure
that the moments correspond to measures with appropriate supports.
For $d$ finite, problem (\ref{lmi}) is a finite-dimensional convex
LP problem in the cone of positive semidefinite matrices, or SDP problem,
for which off-the-shelf solvers are available, in particular implementations
of the primal-dual interior-point methods described in \cite{nesterov94}.

\subsection{Convergence}

The hierarchy of LMI relaxations (\ref{lmi}) generates an asymptotically
converging monotonically increasing sequence of lower bounds on LP (\ref{lp}), i.e.
$J_d \leq J_{d+1}$ for all $d=1,2\ldots$ and $\lim_{d\to\infty} J_d = J_{\infty}$.
Obviously, the relaxed LP cost in problem (\ref{lp})
is a lower bound on the original cost in problem (\ref{ode}), i.e. $J_{\infty} \leq J$.
Under mild assumptions, we can show that actually $J_{\infty} = J$,
but the theoretical background required to prove this lies beyond
the scope of this application paper.

\subsection{Long range behavior}

In the case terminal time $T$ tends to infinity, normalization
constraint (\ref{norm}) becomes irrelevant, and the overall mass of
occupation measure $\mu$ tends to infinity. A more
appropriate formulation consists then of normalizing measure $\mu$
to a probability measure
\[
\pi = \frac{\mu}{T} = \sum_k \frac{\mu_k}{T} = \sum_k \pi_k
\]
so that PDE (\ref{ipde}) becomes
\[
\mathrm{div}(f\pi) = \lim_{T \to \infty} \frac{\mu_0-\mu_T}{T} = 0
\]
and LP problem (\ref{lp}) becomes
\[
\begin{array}{ll}
\inf_{\pi} & \sum_k \int h(t,x) d\pi_k(x) \\
\mathrm{s.t.} & \sum_k \int \frac{\partial v(t,x)}{\partial t} \:d\pi_k(x) \:\:+ \\
& \sum_k \int Dv(t,x) \cdot f_k(t,x) \:d\pi_k(t,x) \:\:= 0 \\
\end{array}
\]
where test functions $v(t,x)$ satisfy appropriate integrability
and/or periodicity properties. Measure $\pi$ is called
an invariant probability measure and it encodes equilibrium points,
periodic orbits, ergodic behavior, possibly chaotic attractors,
see e.g. \cite{lasota85}, \cite{diaconis99}, \cite{hernandez03} and \cite{gaitsgory09}.

\section{Application to orbital launcher control law validation}

In this section we report the application of our method to simplified
models of an attitude control system (ACS) for a launcher in exo-atmospheric phase.
We consider first a one degree-of-freedom (1DOF) model, and then a 3DOF model.
The full benchmark is described in \cite[Section 9]{safev} but it is not
publicly available, and our computer codes cannot be distributed either.

\subsection{ACS 1DOF}

First we consider the simplified 1DOF model of launcher ACS
in orbital phase. The closed-loop system must follow a given
piecewise linear angular velocity profile. The original benchmark includes a time-delay
and a pulsation width modulator in the control loop, but in our incremental approach
we just discard these elements in a first approximation. See below for
a description of possible extensions of our approach to time-delay systems.
As to the pulsation width modulator, it can be handled if appropriately
modeled as a parametric uncertainty.

\subsubsection{Model}

The system is modeled as a double integrator
\[
I \ddot{\theta}(t) = u(t)
\]
where $I$ is a given constant inertia and $u(t)$ is the torque control.
We denote
\[
x(t) = \left[\begin{array}{c} \theta(t) \\ \dot{\theta}(t) \end{array}\right]
\]
and we assume that both angle $x_1(t)$ and angular velocity $x_2(t)$
are measured, and that the torque control is given by
\[
u(x(t)) = \mathrm{sat}(K^T \mathrm{dz}(x_r(t)-x(t)))
\]
where $x_r(t)$ is the reference signal, $K \in {\mathbb R}^2$ is a given state feedback,
$\mathrm{sat}$ is a saturation function such that $\mathrm{sat}(y) = y$ if
$|y|\leq L$ and $\mathrm{sat}(y) = L\:\mathrm{sign}(y)$ otherwise, $\mathrm{dz}$ is
a dead-zone function such that $\mathrm{dz}(x) = 0$ if $|x_i| \leq D_i$ for some $i=1,2$
and $\mathrm{dz}(x) = 1$ otherwise. Thresholds $L>0$, $D_1>0$ and $D_2>0$ are given.

We would like to verify whether the system state $x(t)$
reaches a given subset $X_T = \{x \in {\mathbb R}^2 \: :\: x^Tx \leq \varepsilon\}$
of the deadzone region after a fixed time $T$, and for all possible initial conditions $x(0)$ chosen
in a given subset $X_0$ of the state-space, and for zero reference signals.
We formulate an optimization problem (\ref{ode}) with systems dynamics defined
as locally affine functions in three cells $X_k$, $k=1,2,3$
corresponding respectively to the linear regime of the torque saturation
\[
X_1 = \{x \in {\mathbb R}^2 \: :\: |K^T x| \leq L\}, \quad f_1(x) = \left[\begin{array}{c}
x_1 \\ - K^T x \end{array}\right]
\]
the upper saturation regime
\[
X_2 = \{x \in {\mathbb R}^2 \: :\: K^T x \geq L\}, \quad f_2(x) = \left[\begin{array}{c}
x_1 \\ - L \end{array}\right]
\]
and the lower saturation regime
\[
X_3 = \{x \in {\mathbb R}^2 \: :\: K^T x \leq -L\}, \quad f_3(x) = \left[\begin{array}{c}
x_1 \\ L \end{array}\right].
\]
The objective function has no integral term and a concave quadratic
terminal term $h_T(x) = -x(T)^Tx(T)$ which we would like to minimize,
so as to find trajectories with terminal states of largest norm.
If we can certify that for every initial state $x(0)$ chosen in $X_0$
the final state $x(T)$ belongs to set included in the
deadzone region, we have validated our controlled system.

\subsubsection{Validation script}

The resulting GloptiPoly 3 script,
implementing some elementary scaling
strategies to improve numerical behavior of the SDP solver,
is as follows:

\small
\begin{verbatim}
I = 27500; % inertia
kp = 2475; kd = 19800; % controller gains
L = 380; % input saturation level 
dz1 = 0.2*pi/180; dz2 = 0.05*pi/180; % deadzone levels
thetamax = 50; omegamax = 5; % bounds on initial conditions
epsilon = sqrt(1e-5); % bound on norm of terminal condition
T = 50; % final time

d = input('order of relaxation ='); d = 2*d;

% measures
mpol('x1',2); m1 = meas(x1); % linear regime
mpol('x2',2); m2 = meas(x2); % upper sat
mpol('x3',2); m3 = meas(x3); % lower sat
mpol('x0',2); m0 = meas(x0); % initial
mpol('xT',2); mT = meas(xT); % terminal

% dynamics on normalized time range [0,1]
% saturation input y normalized in [-1,1]
K = -[kp kd]/L;
y1 = K*x1; f1 = T*[x1(2); L*y1/I]; % linear regime
y2 = K*x2; f2 = T*[x2(2); L/I]; % upper sat
y3 = K*x3; f3 = T*[x3(2); -L/I]; % lower set

% test functions for each measure = monomials
g1 = mmon(x1,d); g2 = mmon(x2,d); g3 = mmon(x3,d);
g0 = mmon(x0,d); gT = mmon(xT,d);

% unknown moments of initial measure
y0 = mom(g0);

% unknown moments of terminal measure
yT = mom(gT);

% input LMI moment problem
cost = mom(xT'*xT);
Ay = mom(diff(g1,x1)*f1)+...
     mom(diff(g2,x2)*f2)+...
     mom(diff(g3,x3)*f3); % dynamics
% trajectory constraints
X = [y1^2<=1; y2>=1; y3<=-1]; 
% initial constraints
X0 = [x0(1)^2<=thetamax^2, x0(2)^2<=omegamax^2];
% terminal constraints
XT = [xT'*xT<=epsilon^2];
% bounds on trajectory
B = [x1'*x1<=4; x2'*x2<=4; x3'*x3<=4]; 

% input LMI moment problem
P = msdp(max(cost), ...
         mass(m1)+mass(m2)+mass(m3)==1, ...
         mass(m0)==1, ...
         Ay==yT-y0, ...
         X, X0, XT, B);

% solve LMI moment problem
[status,obj] = msol(P)

\end{verbatim}
\normalsize

\subsubsection{Numerical results}

With this Matlab code and the SDP solver SeDuMi 1.3 (see
{\tt sedumi.ie.lehigh.edu})
we obtain the following sequence of upper bounds (since we
maximize) on the maximum squared Euclidean norm of the final state:
\begin{center}
\begin{tabular}{c|cccc}
relaxation order $d$ & 1 & 2 & 3 & 4 \\\hline
upper bound $J_d$ & $1.0\cdot10^{-5}$ & $1.0\cdot10^{-5}$ & $1.0\cdot10^{-5}$ & $1.0\cdot10^{-5}$ \\
CPU time (sec.) & 0.2 & 0.5 & 0.7 & 0.9 \\
number of moments & 30 & 75 & 140 & 225
\end{tabular}
\end{center}
In the table we also indicate the CPU time (in seconds, on a standard desktop computer)
and the total number of moments (size of vector $y$ in the LMI relaxation (\ref{lmi}).
We see that the bound obtained at the first relaxation ($d=1$)
is not modified for higher relaxations.
This clearly indicates that all initial conditions are
captured in the deadzone region at time $T$, which is
the box $[-2,2]\frac{10^{-1}\pi}{180} \times [-5,5]\frac{10^{-2}\pi}{180}
\supset \{x \in {\mathbb R}^2 \: :\: x^T x \leq 10^{-5}\}$.

\begin{figure}
\begin{center}
\includegraphics[width=\textwidth]{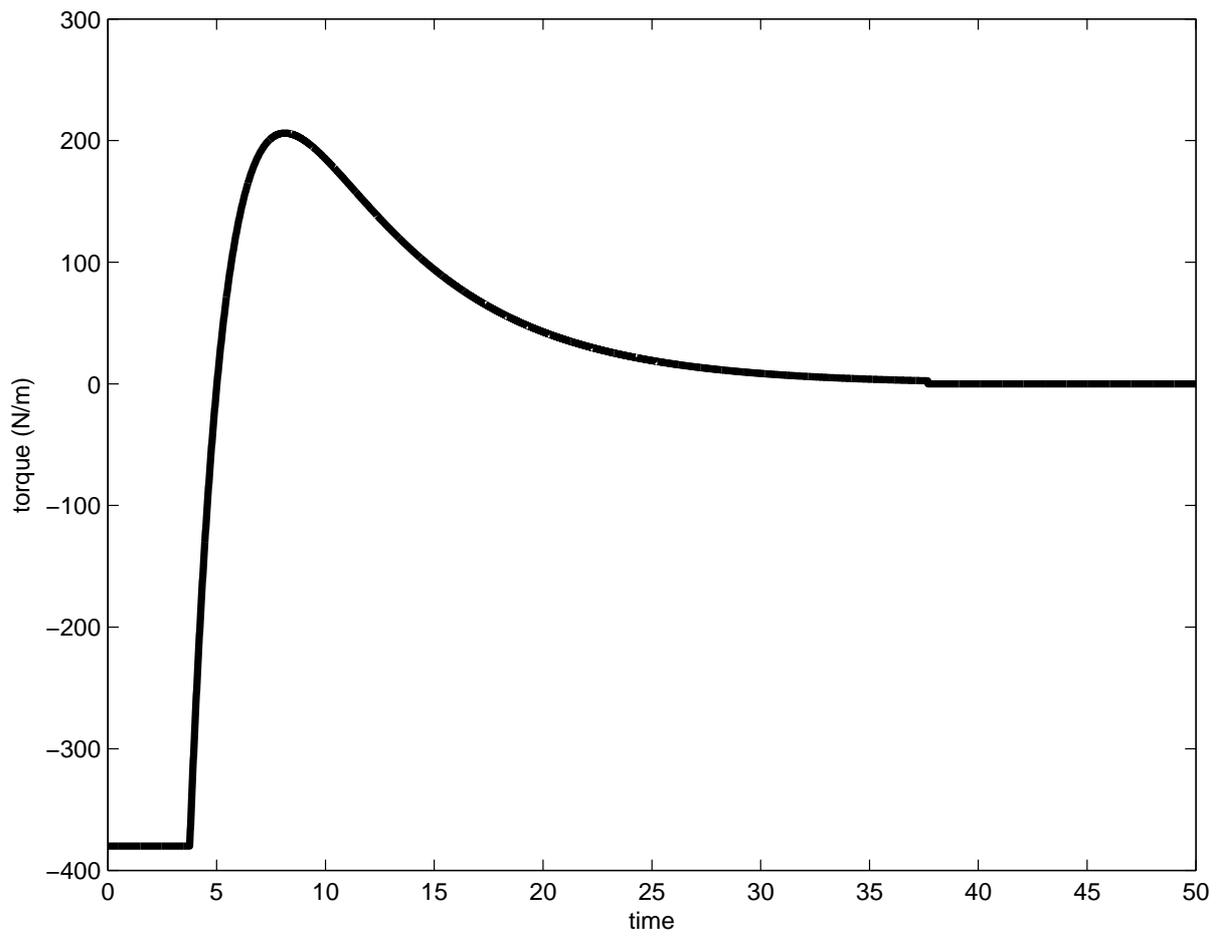}   
\caption{Torque input with lower saturation
during approx. 7\% of the time range.\label{fig}}                              
\end{center}      
\end{figure}

If we want to use this approach to simulate a particular trajectory,
in the code we must modify the definition of the initial measure.
For example for initial conditions $x_1(0) = 50$, $x_2(0) = -1$,
we must insert the following sequence:
\small
\begin{verbatim}
% given moments of initial measure = Dirac at x0
p = genpow(3,d); p = p(:,2:end); % powers
theta0 = 50; omega0 = -1; % in degrees
y0 = ones(size(p,1),1)*[theta0 omega0]*pi/180;
y0 = prod(y0.^p,2);
\end{verbatim}
\normalsize
As previously, the sequence of bounds on the maximum squared Euclidean norm 
of the final state is constantly equal to $1.0\cdot10^{-5}$, and in the
following table we represent as functions of the relaxation order $d$
the masses of measures $\mu_k$, $k=1,2,3$
which are indicators of the time spent by the trajectory in the
respective linear, upper saturation and lower saturation regimes:
\begin{center}
\begin{tabular}{c|cccccccc}
relaxation order $d$ & 1 & 2 & 3 & 4 & 5 & 6 & 7 \\ \hline
$\int d\mu_1$ & 37 & 89 & 92 & 92 & 93 & 93 & 93\\ 
$\int d\mu_2$ & 32 & 5.3 & 0.74 & 0.30 & 0.21 & 0.15 & 0.17 \\
$\int d\mu_3$ & 32 & 5.1 & 7.1 & 6.9 & 6.8 & 6.9 & 7.0
\end{tabular}
\end{center}
This indicates that most of the time (approx. 93\%) is spent in
the linear regime, with approx. 7\% of the time spent in the lower
saturation regime, and a negligible amount of time is spent in the
upper saturation regime. This is confirmed by simulation, see Figure \ref{fig}.

\subsection{Dealing with uncertainty}

Finally, note that we can incorporate real parametric uncertainty
in the dynamics, if required. Each uncertain parameter must be
introduced as an additional state of the system, and this generally
makes the dynamics polynomial in the extended state.
The parameters must be constrained
to an explicitly given compact semialgebraic set, and
the overall trajectory optimization problem consists of
finding the worst-case uncertain parameter instance.

For illustration, we modify the previous script to cope
with uncertainty entering affinely the dynamics. We assume
that the (reciprocal of the) inertia is subject to multiplicative
relative uncertainty, i.e. we replace occurences of $I$ with $\frac{1}{1+u}I$
for $u$ a real parameter such that $u^2 \leq U$ with $U>0$
a given threshold, say $U=\frac{1}{2}$. The resulting
GloptiPoly 3 script is as follows:

\small
\begin{verbatim}

I = 27500; % inertia
kp = 2475; kd = 19800; % controller gains
L = 380; % input saturation level 
dz1 = 0.2*pi/180; dz2 = 0.05*pi/180; % deadzone levels
thetamax = 50; omegamax = 5; % bounds on initial conditions
epsilon = sqrt(1e-5); % bound on norm of terminal condition
T = 50; % final time

ui = 0.5; % uncertainty level (in percent)

d = input('order of relaxation ='); d = 2*d;

% measures = states + uncertain parameter
mpol('x1',2); mpol u1; m1 = meas(x1,u1); % linear regime
mpol('x2',2); mpol u2; m2 = meas(x2,u2); % upper sat
mpol('x3',2); mpol u3; m3 = meas(x3,u3); % lower sat
mpol('x0',2); mpol u0; m0 = meas(x0,u0); % initial
mpol('xT',2); mpol uT; mT = meas(xT,uT); % terminal

% dynamics on normalized time range [0,1]
% saturation input y normalized in [-1,1]
K = -[kp kd]/L;
y1 = K*x1; f1 = T*[x1(2); (1+u1)*L*y1/I]; % linear regime
y2 = K*x2; f2 = T*[x2(2); (1+u2)*L/I]; % upper sat
y3 = K*x3; f3 = T*[x3(2); -(1+u3)*L/I]; % lower set

% test functions for each measure = monomials
g1 = mmon([x1;u1],d); g2 = mmon([x2;u2],d); g3 = mmon([x3;u3],d);
g0 = mmon([x0;u0],d); gT = mmon([xT;uT],d);

% unknown moments of initial measure
y0 = mom(g0);

% unknown moments of terminal measure
yT = mom(gT);

% input LMI moment problem
cost = mom(xT'*xT);
Ay = mom(diff(g1,x1)*f1)+...
     mom(diff(g2,x2)*f2)+...
     mom(diff(g3,x3)*f3); % dynamics
% trajectory constraints
X = [y1^2<=1; y2>=1; y3<=-1]; 
% initial constraints
X0 = [x0(1)^2<=thetamax^2, x0(2)^2<=omegamax^2];
% terminal constraints
XT = [xT'*xT<=epsilon^2];
% bounds on trajectory
B = [x1'*x1<=4; x2'*x2<=4; x3'*x3<=4]; 
% bounds on uncertain parameter
U = [u1^2<=ui^2; u2^2<=ui^2; u3^2<=ui^2];
U0 = [u0^2<=ui^2]; UT = [uT^2<=ui^2];

% input LMI moment problem
P = msdp(max(cost), ...
         mass(m1)+mass(m2)+mass(m3)==1, ...
         mass(m0)==1, ...
         Ay==yT-y0, ...
         X, X0, XT, B, ...
         U, U0, UT);

% solve LMI moment problem
[status,obj] = msol(P)
\end{verbatim}
\normalsize
Introducing an uncertain parameter amounts to adding a state
variable to the model, and this has a significant impact
on the overall computational time, as shown in the table below:
\begin{center}
\begin{tabular}{c|cccc}
relaxation order $d$ & 1 & 2 & 3 & 4 \\\hline
upper bound $J_d$ & $1.0\cdot10^{-5}$ & $1.0\cdot10^{-5}$ & $1.0\cdot10^{-5}$ & $1.0\cdot10^{-5}$ \\
CPU time (sec.) & 0.8 & 2.4 & 7.3 & 17.7 \\
number of moments & 75 & 224 & 501 & 946
\end{tabular}
\end{center}
Note however that for this example there is no effect
on the bounds, and the control law is also validated 
in the presence of uncertainty on inertia.

\subsection{ACS 3DOF}

The main difficulty with the 3 degree-of-freedom version of the ACS
benchmark is the large number of states (4 quaternions for the launcher
attitude, 3 angular velocities, and 2 additional states for the quaternion
reference signal) and the non-linear
(quadratic) dynamics.

\subsubsection{Model}

The following example illustrates the validation of a given
control law following a given constant roll velocity equal to $w^R_1=20^o/s$,
in the absence of actuator saturation, but in the presence
of non-linear coupling between the three axes of the launcher.
System dynamics are given by $\dot{x}=f(x)$ with
\[
x = \left[\begin{array}{c}q\\w\\q^R\end{array}\right], \quad
f(x) = \left[\begin{array}{c}f_q(x)\\f_w(x)\\f_{q^R}(x)\end{array}\right],
\]
and
\[
f_q(x) = \frac{1}{2}Q(x)\left[\begin{array}{c}0\\w\end{array}\right],
\]
\[
f_w(x) = I^{-1}(u(x)-\Omega(x)Iw),
\]
\[
f_{q^R}(x) = \left[\begin{array}{r}
-\frac{1}{2}q^R_1w^R_1 \\
\frac{1}{2}q^R_0w^R_1
\end{array}\right]
\]
where
\[
Q(x) = \left[\begin{array}{rrrr}
q_0 & -q_1 & -q_2 & -q_3 \\
q_1 & q_0 & -q_3 & q_2 \\
q_2 & q_3 & q_0 & -q_1 \\
q_3 & -q_2 & q_1 & q_0
\end{array}\right], \quad
\Omega(x) = \left[\begin{array}{rrr}
0 & -w_3 & w_2 \\ w_3 & 0 & -w_1 \\ -w_2 & w_1 & 0
\end{array}\right],
\]
\[
u(x) = -K^D w^E - K^P q^E,
\]
\[
K^D = \left[\begin{array}{ccc}
K^D_1 & 0 & 0 \\
0 & K^D_2 & 0 \\
0 & 0 & K^D_3 \\
\end{array}\right], \quad
K^P = \left[\begin{array}{ccc}
0 & 0 & 0 \\
0 & K^P_2 & 0 \\
0 & 0 & K^P_3 \\
\end{array}\right],
\]
\[
w^E = \left[\begin{array}{c}
w_1 - w^R_1 \\
w_2 \\
w_3
\end{array}\right], \quad
q^E = 2\left[\begin{array}{cccc}
q^R_0 & q^R_1 & 0 & 0 \\ -q^R_1 & q^R_0 & 0 & 0 \\
0 & 0 & q^R_0 & q^R_1 \\ 0 & 0 & -q^R_1 & q^R_0
\end{array}\right] q,
\]
and $q \in {\mathbb R}^4$ is the quaternion modeling the attitude,
$w \in {\mathbb R}^3$ models the angular velocities,
$q^R \in {\mathbb R}^2$ is the reference quaternion to follow,
which depends on the constant reference velocity $w^R_1 \in {\mathbb R}$,
$q^E \in {\mathbb R}^4$ is the quaternion error, and
$w^E \in {\mathbb R}^3$ is the angular velocity error.
Both vectors $q(t)$ and $q^R(t)$ have constant Euclidean norm
since $\frac{d}{dt}\|q\|^2_2 = 2q^T\dot{q} = q^Tf_q(x) = 0$ for all $x$
and similarly $\frac{d}{dt}\|q^R\|^2_2 = 0$, so we enforce
the algebraic constraints
\[
\|q\|^2_2 = 1, \quad \|q^R\|^2_2 = 1.
\]
The proportional-derivative control law $u(x)$ is designed
so that velocity $w_1(t)$ follows reference velocity $w^R_1$, and
our validation tasks consists of optimizing over the worst-case
trajectory starting at a given initial condition $x(0)$
and maximizing the objective function
\[
J = \int_0^T (w_1(t)-w^R_1)^2 dt
\]
for a given time horizon $T$. If we can guarantee a sufficiently
small upper bound $J_d$ on this objective function, for a given
relaxation order $d=1,2,\ldots$, the control law
is validated.

\subsubsection{Validation script}

We use the following GloptiPoly 3 script (note however
that the function {\tt G2I} to generate the inertia data
is not available for the reader):

\small
\begin{verbatim}
% Inertia
Ig = [ 27500,   -50, -1100;
         -50, 45300,  -220;
       -1100,  -220, 44100];
[Ge2In,Ip] = G2I(Ig);

% controller gains
w0 = 0.3; xi = 0.7;
Kp1 = Ip(1,1)*w0^2; Kd1 = 2*Ip(1,1)*xi*w0;
w0 = 0.25; xi = 3;
Kp2 = Ip(2,2)*w0^2; Kd2 = 2*Ip(2,2)*xi*w0;
Kp3 = Ip(3,3)*w0^2; Kd3 = 2*Ip(3,3)*xi*w0;

% initial conditions
Q0 = [1;0;0;0]; % attitude quaternion
W0 = [0;0;0].*(pi/180); % angular velocities

T = 50; % final time

% reference velocity
WR1 = 20*pi/180;

d = input('order of relaxation =');

% occupation measure
mpol('q',4); mpol('w',3); mpol('qr',2);
x = [q;w;qr]; m = meas(x); 

% initial measure
mpol('q0',4); mpol('w0',3); mpol('qr0',2);
x0 = [q0;w0;qr0]; m0 = meas(x0); 

% terminal measure 
mpol('qt',4); mpol('wt',3); mpol('qrt',2);
xt = [qt;wt;qrt]; mt = meas(xt);
 
% dynamics on normalized time range [0,1]
Omega = [0 -w(3) w(2); w(3) 0 -w(1); -w(2) w(1) 0];
Qr = [qr(1) qr(2) 0 0;
      -qr(2) qr(1) 0 0;
      0 0 qr(1) qr(2);
      0 0 -qr(2) qr(1)];
dq = 2*Qr*q;
Q = [q(1) -q(2) -q(3) -q(4);
     q(2) q(1) -q(4) q(3);
     q(3) q(4) q(1) -q(2);
     q(4) -q(3) q(2) q(1)];
fq = 1/2*Q*[0;w];
Torque = [-Kd1*(w(1)-WR1);
          -Kd2*(w(2)-0)-Kp2*2*dq(3);
          -Kd3*(w(3)-0)-Kp3*2*dq(4)];
fw = Ig\(Torque-Omega*Ig*w);
fqr = [-1/2*qr(2)*WR1;1/2*qr(1)*WR1];
f = T*[fq;fw;fqr];

% test functions
g = mmon(x,d);
g0 = mmon(x0,d);
gt = mmon(xt,d);

% given moments of initial measure (Dirac)
p = genpow(10,d); p = p(:,2:end); % powers
y0 = ones(size(p,1),1)*[Q0' W0' 1 0];
y0 = prod(y0.^p,2);

% unknown moments of final measure
p = genpow(10,d); p = p(:,2:end); % powers
yt = ones(size(p,1),1)*[qt' wt' qrt'];
yt = mom(prod((yt.^p)')');

% input LMI moment problem
cost = mom((w(1)-WR1)^2); % objective function
Ay = mom(diff(g,x)*f); % dynamics
X = [q'*q == 1; qr'*qr == 1]; % quaternion
Xt = [qt'*qt == 1; qrt'*qrt == 1]; % final quaternion
P = msdp(max(cost), mass(m)==1, yt==Ay+y0, X, Xt);

% solve LMI moment problem
[status,obj] = msol(P)
\end{verbatim}
\normalsize

\subsubsection{Numerical results}

We obtain the monotonically decreasing
upper bounds $J_d$, $d=1,\ldots,4$ reported in the following table:
\begin{center}
\begin{tabular}{c|cccc}
relaxation order $d$ & 1 & 2 & 3 & 4 \\\hline
upper bound $J_d$ & $\infty$ & $7.5198\cdot 10^{-3}$ & $2.9110\cdot 10^{-3}$ & $2.9085\cdot 10^{-3}$ \\
CPU time (sec.) & 0.2 & 11.6 & 24.2 & 2640 \\
number of moments & 110 & 770 & 1430 & 5720
\end{tabular}
\end{center}
The first LMI relaxation (110 moments of degree up to 2)
seems to be unbounded above (the dual LMI problem is infeasible) and 
hence it conveys no useful information.
The second LMI relaxation (770 moments of degree up to 3)
and the third LMI relaxation (1430 moments of degree up to 4)
are solved with SeDuMi 1.3 in a few seconds.
The fourth LMI relaxation is more challenging, with 5720 moments
of degree up to 5, and it requires less than one hour of CPU time
on a standard PC. We observe however that a useful upper bound
on $J$ is already obtained at a low relaxation order (say $d=2$
or $d=3$) at a low computational cost, and
that the computational burden increases significantly for $d=4$
without dramatic improvements in the quality of the upper bound.

\subsubsection{Coping with saturations and dead-zones}

We can modify the above code to cope with actuator saturation
and dead-zone, proceeding exactly as described above for the ACS
1DOF benchmark problem, namely by splitting the trajectory occupation measures
into local occupation measures defined on semialgebraic
(here polyhedral) cells. In the presence of saturations
on the three actuators, this generates $3^3=27$ local
occupation measures. This should not be an issue from
the computational point of view, the limiting factor
being essentially the size of the largest SDP block,
which grows polynomially with the relaxation order,
but with an exponent which is the number of variables
(here equal to 10).

\subsubsection{Following a time-varying reference signal}

Finally, if we want to validate a control law to follow
a time-varying velocity profile (instead of a constant velocity
as above), we must introduce time as an additional variable
in the trajectory occupation measure, we must introduce
the velocity as an additional state, and we must split the
occupation measure into local occupation measures, each
of which corresponding to given time-invariant dynamics.
Alternatively, we can also model time-varying dynamics,
but the dependence on time must be polynomial (which is
not the case e.g. if the reference signal is piecewise
linear).

\section{Discussion}\label{discussion}

In this section we discuss the limitations of our
approach and we also describe possible extensions.

\subsection{Limitations}

The main limiting factor for our method is the total number of variables
(number of states plus number of uncertain parameters), since the
computational burden of solving the LMI relaxations
grows polynomially as a function of the LMI relaxation order, but the order
of the polynomial dependence depends linearly on the number of variables
in the measures (here the number of states). So the critical dimension
is not too much the number of measures (that is the number of cells partioning
the state-space, which corresponds to our models of nonlinearities),
or the degree of the polynomials in the cost and/or dynamics,
but rather the number of states. Systems with a high number of states can
be handled with these techniques, but only when exploiting problem structure
and sparsity.

More explicitly, a rough complexity analysis can be carried out as follows.
If an LMI relaxation has the simplified form
\[
\begin{array}{ll}
\inf_y & c^T y \\
\mathrm{s.t.} & M_d(y) \succeq 0
\end{array}
\]
where $M_d(y)$ is the moment matrix of a measure of $n$ variables at relaxation order
$d$, then the number of variables in vector $y$ is
\[
N=\left(\begin{array}{c}n+2d\\n\end{array}\right) = \frac{(n+2d)!}{n!\,(2d)!}
\]
and the size of matrix $M_d(y)$ is
\[
M=\left(\begin{array}{c}n+d\\n\end{array}\right) = \frac{(n+d)!}{n!\,d!}.
\]
A standard primal-dual interior-point algorithm to solve
this LMI at given relative accuracy $\epsilon>0$ (duality gap threshold)
requires a number of iterations (Newton steps) growing as $O(M^{\frac{1}{2}}\log\epsilon)$,
whose dependence on $M$ is sublinear, hence almost negligible. In practice, on
well-conditioned problems, we observe that the number of Newton iterations
is between 5 and 50, see \cite{nesterov94}, \cite{boyd} and \cite{bental}.

In the real model of computation (for which each addition, subtraction, multiplication,
division of real numbers has unit cost), each Newton iteration requires $O(N^2M)+O(N^3M)+O(N^2M^2)$ operations
to form the linear system of equations, and $O(M^3)$ operations 
to solve the system and find the search direction. When solving a hierarchy
of simple LMI relaxations as described above, the number of variables $n$ is fixed,
and the relaxation order $d$ varies, so the dominating term in the complexity
estimate grows in $O(d^{4n})$, which clearly shows a strong dependence
on the number of variables. Even though the growth of the computational burden
is polynomial in the relaxation order, the exponent is 4 times the number
of variables. One should however keep in mind that these estimates
are (usually very loose) asymptotic upper bounds, and that the observed
computational complexity grows much more moderately in practice (at least
in the absence of conditioning and numerical stability issues).

Finally, if the LMI relaxation has the form
\[
\begin{array}{ll}
\inf_y & c^T y \\
\mathrm{s.t.} & M_d(y_k) \succeq 0, \quad k=1,2,\ldots,K
\end{array}
\]
where $M_d(y_k)$ is the moment matrix of a measure $\mu_k$ of $n$ variables at relaxation $d$,
and we have $K$ measures, then the above complexity estimate grows in $O(Kd^{4n})$.
The impact on the computation burden of the number $n$ of variables in each measure $\mu_k$
is thus much more critical than the number $K$ of measures. In our target application,
$n$ is the number of states and uncertain parameters, $K$ is the number
of cells used to model nonlinearities, and a lower bound on the minimum relaxation order $d$
is given by the degree of the polynomial data (dynamics, constraints, objective function).
 
\subsubsection{Accuracy}

The original validation problem can be formulated as infinite-dimensional
linear programming problems on the dual cones of nonnegative measures
and continuous functions, so that there is no hope of finite convergence
of finite-dimensional LMI optimization techniques. Convergence is
guaranteed only asymptotically. Obviously, the accuracy depends
on the speed of convergence of the hierarchy of LMI problems, but
this speed is impossible to evaluate a priori. The only guarantee
is that we have a monotically increasing sequence of lower bounds
on the objective function to be minimized.

As far as solving finite-dimensional LMI problems is concerned,
we should emphasize the fact that there is currently no proof
of backward stability of implementations of LMI solvers. Moreover,
evaluation of the conditioning of a given LMI problem is difficult,
and estimates can be obtained only at the price of solving several
instances of (slightly modified versions) the original LMI problem.
To certify the output of a numerical algorithm, we must simultaneously
ensure backward stability of the algorithm and well-conditioning of
the problem. None of these two properties can be ensured when solving LMI
problems, given the current state-of-the-art in LMI solvers.
But this is not specific to our method, and any validation-verification
technique which is based on LMI optimization is subject to the
same limitations.

\subsubsection{Time-delay systems}

Our techniques can in principle cope with time-delays, but this requires a
significant modification of the approach described in this document. Let us only
sketch the main ideas, in the case of an ordinary differential equation with
one time-delay
\[
\dot{x}(t) = f(x(t)) + g(x(t-\tau)), \quad \forall t \in [0,T]
\]
with boundary conditions
\[
\quad x(t) = \xi(t), \:\:\forall t \in [-\tau,0]
\]
where $\xi(t)$ is a given function recording the state history due to the given delay
$\tau \in \mathbb R$.
Instead of transporting a probability measure $\mu_t(dx)$ supported on $X \subset {\mathbb R}^n$
from initial time $t=0$ to terminal time $t=T$, we must transport the state history
in an occupation measure $\mu_t(ds,dx)$ supported on $[-\tau,0]\times X$
for $t \in [0,T]$.
 
\subsubsection{Discrete-time systems}

In this paper we deal only with continuous-time
systems. Discrete-time systems of the form
\begin{equation}\label{discrete}
x_{k+1} = f(x_k)
\end{equation}
must be handled differently. Denoting by $\mu_k(x)$ the probability measure
transported along dynamical system (\ref{discrete}), the discrete-time analogue of
Liouville's transport equation (\ref{pde}) reads
\begin{equation}\label{dmeas}
\mu_{k+1}({\mathcal X}) = \int_{f^{-1}({\mathcal X})}\mu_k(dx) = \int I_{\mathcal X}(f(x))\mu_k(dx)
\end{equation}
where $\mathcal X$ is any subset of the sigma-algebra of state set $X$,
and $I_{\mathcal X}$ is the indicator function equal to one in $\mathcal X$
and zero outside. It follows from (\ref{dmeas})
that the moments of measure $\mu_{k+1}$ can
be expressed linearly as functions of moments of measure $\mu_k$.
Besides this analogy, the resulting discrete-time generalized moment problem
differs significantly from its continuous-time counterpart. For this reason, handling
discrete-time systems would require an important modification of the approach described
in this document.

\section{Conclusion}

In this paper we do not follow the mainstream Lyapunov approach to dynamical
systems stability and/or performance validation. \cite{lyapunov92}
originally introduced his approach to conclude about systems behaviour
while avoiding explicit computations of system trajectories. In the 1990s
it was combined with LMI techniques, see e.g. \cite{boyd94,nesterov94} 
to provide numerically quadratic certificates of stability and/or performance
for linear and nonlinear systems. Later on, it was extended to more general,
but still numerical polynomial sum-of-squares (SOS) certificates, see e.g. \cite{henrion05}.
By contrast, our approach is genuinely primal, in the sense that we
directly optimize numerically over systems trajectories, using measures,
moments and LMI techniques. Indeed, since the computation of SOS certificates
is numerical anyway, we believe that in the context of validation,
LMI techniques should rather be used
to optimize systems trajectories, instead of optimizing indirect certificates
for these trajectories. Since we are using primal-dual SDP solvers to solve
the hierarchy of LMI problems, solving the dual to the problem of optimizing trajectories
provides anyway certificates of feasibility or infeasibility, i.e. of stability,
instability and/or performance.

As briefly described in \S \ref{discussion}, the measure/LMI approach can be extended
easily to systems with real uncertainties, as soon as the uncertainties enter
polynomially in the dynamics, and they are bounded in explicitly given
semialgebraic sets (e.g. balls or boxes). See e.g. \cite[Section 13.1]{lasserre09}
for more information on robust optimization. As usual, the number of real uncertain
parameters should be kept small enough to ensure computational tractability.

Note that in this paper the set of initial conditions is assumed to be given,
and when the initial measure is unknown, its support is constrained to be included
in the set of initial conditions. In the context of validation/verification it
could be relevant to maximize the size of the set of initial conditions,
and we are currently investigating this problem from the point of view of
occupation measures. 

The measure/LMI approach can readily deal with systems with piecewise polynomial
(or rational) models, e.g. systems with input/output saturations and dead-zones.
It can be seen as a (primal) extension of early attempts of use of LMI/Lyapunov
techniques in the context of systems with saturations, see e.g. \cite{henrion99}
and more recently \cite{garulli11} and \cite{tarbouriech11}.
As briefly discussed in \S \ref{discussion}, our approach can be extended without major difficulty
to time-delay, stochastic, discrete-time systems or
hybrid system dynamics where transitions between models are ruled e.g. by
probability measures, see e.g. \cite{diaconis99}, \cite{hernandez03}
or \cite{barkley06}.

\section*{Acknowledgment}

We are grateful to Dimitri Peaucelle and Prathyush Menon
for their comments on this work.

\end{document}